\documentclass[12pt,french]{smfart}
\usepackage{amsthm}

\usepackage[T1]{fontenc}
\usepackage[english,francais]{babel}

\begin{document}

\title{Courte preuve d'une identit\'{e} comportant les nombres de Bernoulli}
\date{\today}
\author{Abdelmoum\`ene ZEKIRI}
\address{Facult\'e de Math\'emath\'ematiques, U.S.T.H.B., LA3C \\
Bp 32 El Alia 16111 Bab Ezzouar Alger, Alg\'erie.}
\email{azekiri@usthb.dz,czekiri@gmail.com}
\author{Farid BENCHERIF}
\address{Facult\'e de Math\'emath\'ematiques, U.S.T.H.B., LA3C \\
Bp 32 El Alia 16111 Bab Ezzouar Alger, Alg\'erie.}
\email{fbencherif@usthb.dz, fbencherif@yahoo.fr}

\begin{abstract}
Dans cet article, nous donnons une courte preuve d'une relation g\'{e}n\'{e}%
ralisant de nombreuses identit\'{e}s pour les nombres de Bernoulli. 
\end{abstract}
\maketitle

\section{Introduction et r\'{e}sultat}

Soit $S_{x}(z):=\frac{ze^{zx}}{e^{z}-1}$ la s\'{e}rie g\'{e}n\'{e}ratrice
exponentielle des polyn\^{o}mes $(B_{n}(x))_{n\geq 0}$ de Bernoulli. La
suite $(B_{n})_{n\geq 0}$ des nombres de Bernoulli est d\'{e}finie par%
\begin{equation*}
B_{n}=B_{n}(0),\text{ \ }(n\geq 0).
\end{equation*}%
On a alors 
\begin{equation*}
B_{n}(x)=\sum_{k=0}^{n}\binom{n}{k}B_{k}x^{n-k}.
\end{equation*}%
De nombreux auteurs se sont int\'{e}ress\'{e}s \`{a} la remarquable identit%
\'{e} suivante 
\begin{equation}
(-1)^{m}\sum_{k=0}^{m}\binom{m}{k}B_{n+k}-(-1)^{n}\sum_{k=0}^{n}\binom{n}{k}%
B_{m+k}=0.  \label{r1}
\end{equation}%
Cette identit\'{e} pos\'{e}e comme probl\`{e}me par Carlitz \cite{car1} en
1971 fut prouv\'{e}e en 1972 par Shanon \cite{sha} puis par Gessel \cite{ges}
en 2003, par Wu, Sun et Pan \cite{wu} en 2004, par Vassilev-Missana \cite%
{vas} en 2005 et par Chen et Sun \cite{che} en 2009. Gould et J. Quaintance 
\cite{gou} la prouvent de nouveau en 2014. Enfin Prodinger \cite{pro} en
donne une preuve remarquablement tr\`{e}s courte en 2014 par l'emploi d'une s%
\'{e}rie formelle \`{a} deux variables. Il est imm\'{e}diat de constater que
(\ref{r1}) est un cas particulier, obtenue pour $q=0$, de la relation
suivante \'{e}tablie en 2012 par le second auteur et Garici \cite{ben} 
\begin{equation}
(-1)^{m}\sum_{k=0}^{m+q}\binom{m+q}{k}\binom{n+q+k}{q}B_{n+k}-(-1)^{n+q}%
\sum_{k=0}^{n+q}\binom{n+q}{k}\binom{m+q+k}{q}B_{m+k}=0.  \label{r2}
\end{equation}

Dans cette note, nous pr\'{e}sentons une tr\`{e}s courte preuve de la
relation (\ref{r2}) qui repose essentiellement sur le fait bien connu que
tout polyn\^{o}me de Bernoulli $B_{n}(x)$ d'indice $n$ impair s'annule pour $%
x=\frac{1}{2}$ et qu'on prouve facilement en remarquant que $S_{\frac{1}{2}%
}(-z)=S_{\frac{1}{2}}(z)$. En consid\'{e}rant alors la forme lin\'{e}aire $L$
d\'{e}finie sur le $%
\mathbb{Q}
$- espace vectoriel $%
\mathbb{Q}
\left[ x\right] $ par $L(x^{n})=B_{n}$, pour tout $n\geq 0$, cette propri%
\'{e}t\'{e} se traduit par 
\begin{equation*}
L\left( \left( x+\frac{1}{2}\right) ^{2n+1}\right) =B_{2n+1}\left(\frac{1}{2}\right)=0\qquad(n\geq 0).
\end{equation*}%
Il suffit alors de remarquer que le polyn\^{o}me%
\begin{equation*}
P(x)=(-1)^{m+q}x^{n+q}(1+x)^{m+q}-(-1)^{n}x^{m+q}(1+x)^{n+q},
\end{equation*}%
v\'{e}rifie la relation 
\[P(-\frac{1}{2}+x)+(-1)^{q}P(-\frac{1}{2}-x)=0.
\]
 Il en r\'{e}sulte alors que 
 \[
 P^{(q)}(-\frac{1}{2}+x)+P^{(q)}(-\frac{1}{2}-x)=0.
 \]
On en d\'{e}duit que $P^{(q)}(x)$ est un polyn\^{o}me impair en $x+\frac{1}{2}$ et par cons\'{e}quent que 
\[
L\left( P^{(q)}(x)\right) =0. 
\]
Ce qui fournit la relation (\ref{r2}) en constatant que 
\begin{equation*}
\frac{1}{q!}P^{(q)}(x)=(-1)^{m}\sum_{k=0}^{m}\binom{m+q}{k}\binom{n+q+k}{q}x^{n+k}-(-1)^{n+q}\sum_{k=0}^{n}\binom{n+q}{k}\binom{m+q+k}{q}x^{m+k}.
\end{equation*}%
Signalons que la relation (\ref{r2}) g\'{e}n\'{e}ralise de nombreuses identit\'{e}s connues pour les nombres de Bernoulli.

Pour $m=n$ et $q=1$, la relation permet d'obtenir une relation prouv\'{e}e par Ettingshausen \cite{ett} en 1827, red\'{e}couverte par Seidel \cite{sei} en 1877 puis de nouveau par Kaneko \cite{kan} en 1995 et prouv\'{e}e aussi par Gessel \cite{ges} en 2003 ainsi que par Chen et Sun \cite{che} et par Cigler \cite{cig} en 2009.

Pour $q=1$, on obtient la relation \'{e}tablie par Momiyama \cite{mom} en 2001, puis prouv\'{e}e de nouveau par Wu Sun et Pan \cite{wu} en 2004 et par Chen et sun \cite{che} en 2009.

Pour $m=n$ et $q=3$, on retrouve une relation \'{e}tablie en par Chen et Sun \cite{che} en 2009.

Pour $m=n$ et $q$ impair, on retrouve la relation \'{e}tablie par les auteurs \cite{zek} en 2011.

\end{document}